\documentclass[reqno]{amsart}
\usepackage{amsmath,amsthm}

\newtheorem{Theorem}{Theorem}
\newtheorem{Definition}[Theorem]{Definition}

\newtheorem{Proposition}[Theorem]{Proposition}
\newtheorem{Remark}{Remark}

\begin{document}
\title[Lifespan for the periodic DNLS]%
{Lifespan of strong solutions to the periodic derivative nonlinear Schr\"odinger equation}

\author[K. Fujiwara]{Kazumasa Fujiwara}

\address{%
Department of Pure and Applied Physics \\ Waseda University \\
3-4-1, Okubo, Shinjuku-ku, Tokyo 169-8555 \\ Japan}

\email{k-fujiwara@asagi.waseda.jp}

\thanks{
The first author was partly supported by Grant-in-Aid for JSPS Fellows no 16J30008.
}

 \author[T. Ozawa]{Tohru Ozawa}

\address{%
Department of Applied Physics \\ Waseda University \\
3-4-1, Okubo, Shinjuku-ku, Tokyo 169-8555 \\ Japan}

\email{txozawa@waseda.jp}

\maketitle

\renewcommand{\thefootnote}{\fnsymbol{footnote}}

\begin{abstract}
An explicit lifespan estimate is presented
for the derivative Schr\"odinger equations
with periodic boundary condition.
\end{abstract}

\section{Introduction}
We consider the Cauchy problem for the following
derivative nonlinear Schr\"odinger (DNLS) equation:
	\begin{align}
	\begin{cases}
	i \partial_t u + \partial_x^2 u = \lambda \partial_x (|u|^{p-1} u),
	&\qquad t \in \lbrack 0,T), \quad x \in \mathbb T,\\
	u(0) = u_0,
	&\qquad x \in \mathbb T
	\end{cases}
	\label{eq:1}
	\end{align}
on one-dimensional torus $\mathbb T = \mathbb R / 2 \pi \mathbb Z$,
where $p > 1$ and $\lambda \in \mathbb{C} \backslash \{0\}$.
The aim of this paper is to study an explicit upper bound
of lifespan of solutions for \eqref{eq:1}
in terms of the data $u_0$ in the case $\mathrm{Re} \lambda \neq 0$.

The original DNLS equation on $\mathbb R$ with $p=3$ and $\lambda = -i$
with additional terms
was derived in plasma physics for a model of Alfv\'en wave
(see \cite{bib:13,bib:17}).
By a simple computation,
if $\lambda \in i \mathbb R$,
then we have the charge ($L^2$) conservation law
for solutions of \eqref{eq:1} and
	\begin{align}
	\begin{cases}
	i \partial_t u + \partial_x^2 u = \lambda \partial_x (|u|^{p-1} u),
	&\qquad t \in \lbrack 0,T), \quad x \in \mathbb R,\\
	u(0) = u_0,
	&\qquad x \in \mathbb R
	\end{cases}
	\label{eq:2}
	\end{align}
with any $p > 1$.
For the solution $u$ for \eqref{eq:2} with $p=3$ and $\lambda = -i$,
the gauge transformed solution $v$ defined by
	\[
	v(t,x) = u(t,x) \exp \bigg( \frac{i}{2} \int_{-\infty}^x |u(t,y)|^2 dy \bigg)
	\]
satisfies
	\begin{align}
	i \partial_t v + \partial_x^2 v = -i |v|^{2} \partial_x v,
	\qquad t \in \lbrack 0,T), \quad x \in \mathbb R.
	\label{eq:3}
	\end{align}
Similarly, in the case of \eqref{eq:1} with $p=3$ and $\lambda = -i$,
the gauge transformed solution $w$ defined by
	\[
	w(t,x) = u(t,x) \exp \bigg( \frac i 2 \int_{0}^x |u(t,y)|^2 dy
	- \frac i 2 \int_0^t \mathrm{Im} \Big( \overline{u(t',0)} \partial_x u(t',0) \Big)
	+ 4 |u(t',0)|^4 dt' \bigg)
	\]
satisfies
	\begin{align}
	i \partial_t w + \partial_x^2 w = -i |w|^{2} \partial_x w,
	\qquad t \in \lbrack 0,T), \quad x \in \mathbb T.
	\label{eq:4}
	\end{align}
Then the following energies are conserved:
	\begin{align*}
	E_1(v(t))
	&= \int_{\mathbb R} \bigg( |\partial_x v(t)|^2
	+ \frac 1 2 \mathrm{Im} |v(t)|^2 \overline{v(t)} \partial_x v(t) \bigg) dx
	= E_1(v(0)),\\
	E_2(w(t))
	&= \int_{\mathbb T} \bigg( |\partial_x w(t)|^2
	+ \frac 1 2 \mathrm{Im} |w(t)|^2 \overline{w(t)} \partial_x w(t) \bigg) dx
	= E_2(w(0)).
	\end{align*}

The well-posedness of \eqref{eq:2} with $p=3$ and $\lambda = -i$
has been studied by many authors.
For example, Tsutsumi and Fukuda showed local well-posedness in the Sobolev space
$H^s(\mathbb R) = (1-\Delta)^{-s/2} L^2(\mathbb R)$ with $s >3/2$ in \cite{bib:22}.
Moreover, by using the gauge transformation,
the $H^s(\mathbb R)$ local and global well-posedness with $s \geq 1/2$
has been studied in \cite{bib:3,bib:8,bib:9,bib:10,bib:18}.
Furthermore,
Biagioni and Linares showed the $H^s(\mathbb R)$ ill-posedness
with $s < 1/2$ in \cite{bib:2}.
This means $H^{1/2}(\mathbb R)$ gives the sharp criteria
for the local well-posedness for \eqref{eq:3}.
We also refer the reader to \cite{bib:7,bib:12,bib:16} for generalized results.

On the other hand,
the well-posedness of the Cauchy problem \eqref{eq:1}
with $p=3$ and $\lambda \in i \mathbb R$ has also been studied.
Tsutsumi and Fukuda showed local well-posedness in the Sobolev space
$H^s(\mathbb T) = (1-\Delta)^{-s/2} L^2(\mathbb T)$ with $s >3/2$ in \cite{bib:22}
as well as for \eqref{eq:2}.
In \cite{bib:11}, Herr showed the local and global well-posedness
in $H^s(\mathbb T)$ with $s \geq 1/2$ by using the modified gauge transformation.
We also refer the reader to \cite{bib:1,bib:6,bib:14,bib:19,bib:21,bib:23}
for generalized results.

Even though, local and global well-posedness for DNLS equation has been studied,
the blowup of solution for DNLS is still open in a general setting,
where the conservation low is insufficient or fails.
Partial results have been obtained in \cite{bib:20}.
In this article, we study the finite time blowup of solutions for \eqref{eq:1}
by using a simple ODE argument.
See also \cite{bib:4,bib:5,bib:15}.

An obvious global solution for \eqref{eq:1} is $u(t,x) = C$ for $C \in \mathbb C$.
So it is necessary to consider a set of initial data without constants
in order to show the finite time blowup of \eqref{eq:1}.
Here we consider the initial data and solutions
with vanishing total density defined as follows:

\begin{Definition}
For $u_0 \in H^2(\mathbb T)$
satisfying $\int_{\mathbb T} u_0(x) dx = 0$,
$u$ is called a strong solution with vanishing total density
of the Cauchy problem \eqref{eq:1}
if there exists $T \in (0,\infty]$ such that
$u \in C^1([0,T); H^2(\mathbb T))$ satisfies \eqref{eq:1}
and $\int_{\mathbb T} u(t,x) dx = 0$ for any $t \in \lbrack 0,T)$.
\end{Definition}

\begin{Remark}
Formally,
	\[
	\frac{d}{dt}
	\int_{\mathbb T} u(t,x) dx
	= (2 \pi)^{1/2} \frac{d}{dt} \hat u(t,0)
	= -i (2 \pi)^{1/2} \mathfrak F
	\big\lbrack - \partial_x^2 u + \lambda \partial_x (|u|^{p-1} u) \big\rbrack (0)
	= 0.
	\]
This implies that if $\int_{\mathbb T} u_0(x) dx = 0$,
then $\int_{\mathbb T} u(t,x) dx = 0$ for any $t \in \lbrack 0, T)$.
\end{Remark}

In this article,
for $H^2(\mathbb T)$ initial data with vanishing total density,
we assume the existence of strong solutions with vanishing total density.
We define the lifespan $T_0$ of a strong solution $u$ to the Cauchy problem \eqref{eq:1}
by
	\[
	T_0 = \sup \{T > 0; \ u \mbox{ is a strong solution for \eqref{eq:1}} \}.
	\]
Then,
from the ordinary differential inequality for
	$\displaystyle
	\int_{0}^{2 \pi} \int_{0}^{x} u(\cdot,x)
	\overline{u(\cdot,y)} dy \thinspace dx,
	$
we may obtain the following equivalent conditions
for the finite time blowup for \eqref{eq:1}
and estimate of lifespan.

\begin{Proposition}
\label{Proposition:2}
Let $u_0 \in L^2(\mathbb T)$ satisfy $\int_{\mathbb T} u_0(x) dx = 0$.
Then the following statements are equivalent:
\begin{itemize}
\item[(i)]
$u_0$ satisfies
	\begin{align}
	\mathrm{Re} \lambda \cdot
	\mathrm{Im}
	\int_{0}^{2 \pi} u_0 (x) \int_0^x \overline{u_0(y)} dy \thinspace dx
	> 0.
	\label{eq:5}
	\end{align}
\item[(ii)]
There exists $\alpha \in \mathbb C$ such that
	\begin{align}
	\mathrm{Re} \thinspace \alpha \cdot \mathrm{Re} \lambda,
	\quad
	\mathrm{Im} \bigg( \alpha
	\int_{0}^{2 \pi} u_0 (x) \int_0^x \overline{u_0(y)} dy \thinspace dx \bigg)
	> 0.
	\label{eq:6}
	\end{align}
\end{itemize}
If $u_0$ satisfies one of the equivalent conditions above and $u_0 \in H^2(\mathbb T)$,
then the corresponding strong solution with vanishing total density
of the Cauchy problem \eqref{eq:1} blows up in finite time.
Moreover,
the associated lifespan is estimated by
	\[
	T_0
	\leq
	\frac{(2 \pi)^{p}}{(p-1) | \mathrm{Re} \thinspace \lambda|}
	\bigg| \int_0^{2 \pi} u_0(x) \int_0^x \overline{u_0(y)} dy \thinspace dx
	\bigg|^{-\frac{p-1}{2}}.
	\]
\end{Proposition}

\begin{Remark}
For $f \in L^2(\mathbb T)$ with vanishing total density,
	\begin{align*}
	\mathrm{Re} \int_0^{2 \pi} f(x) \int_0^x \overline{f(y)} dy \thinspace dx
	&= \frac{1}{2} \int_0^{2 \pi} \frac{d}{dx} \bigg| \int_0^x f(y) dy \bigg|^2 dx\\
	&= \frac{1}{2} \bigg| \int_0^{2 \pi} f(x) dx \bigg|^2
	= 0.
	\end{align*}
This means
	\[
	\int_0^{2 \pi} f(x) \int_0^x \overline{f(y)} dy \thinspace dx
	\in i \mathbb R
	\]
and implies the equivalence between \eqref{eq:5} and \eqref{eq:6}.
\end{Remark}

\section{proof of proposition \ref{Proposition:2}}
Let
	$\displaystyle
	M(t)
	= \mathrm{Im}
	\bigg( \alpha \int_{0}^{2 \pi} u(t,x)
	\int_{0}^{x} \overline{u(t,y)} dy \thinspace dx \bigg),
	$
where $\alpha$ satisfies \eqref{eq:6}.
Then $M(t) > 0$ for sufficiently small $t$.
By a direct calculation, we have
	\begin{align*}
	\frac{d}{dt} M(t)
	&= \mathrm{Im}
	\bigg( \alpha \int_{0}^{2 \pi} \partial_t u(t,x)
	\int_{0}^{x} \overline{u(t,y)} dy \thinspace dx \bigg)\\
	&+ \mathrm{Im}
	\bigg( \alpha \int_{0}^{2 \pi} u(t,x)
	\int_{0}^{x} \overline{\partial_t u(t,y)} dy \thinspace dx \bigg)\\
	&= I_1 + I_2.
	\end{align*}
By the vanishing total density,
$I_1$ and $I_2$ may be computed as follows:
	\begin{align*}
	I_1 &= - \mathrm{Re}
	\bigg( \alpha \int_{0}^{2 \pi} i \partial_t u(t,x)
	\int_{0}^{x} \overline{u(t,y)} dy \thinspace dx \bigg)\\
	&= - \mathrm{Re}
	\bigg( \alpha \int_{0}^{2 \pi}
	\partial_x (-\partial_x u(t,x) + \lambda (|u(t,x)|^{p-1} u(t,x)))
	\int_{0}^{x} \overline{u(t,y)} dy \thinspace dx \bigg)\\
	&= - \mathrm{Re}
	\bigg( \alpha (-\partial_x u(t,2 \pi) + \lambda (|u(t,2 \pi)|^{p-1} u(t,2 \pi)))
	\int_{0}^{2 \pi} \overline{u(t,y)} dy \bigg)\\
	&+ \mathrm{Re}
	\bigg( \alpha \int_{0}^{2 \pi}
	- \overline{u(t,x)} \partial_x u(t,x) + \lambda |u(t,x)|^{p+1} dx \bigg)\\
	&= \mathrm{Re}
	\bigg( \alpha \int_{0}^{2 \pi}
	- \overline{u(t,x)} \partial_x u(t,x) + \lambda |u(t,x)|^{p+1} dx \bigg),\\
	I_2&= \mathrm{Re}
	\bigg( \alpha \int_{0}^{2 \pi} u(t,x)
	\int_{0}^{x} \overline{i \partial_t u(t,y)} dy \thinspace dx \bigg)\\
	&= \mathrm{Re}
	\bigg( \alpha \int_{0}^{2 \pi} u(t,x)
	\overline{(-\partial_x u(t,x) + \lambda (|u(t,x)|^{p-1} u(t,x)))} dx \bigg)\\
	&- \mathrm{Re}
	\bigg( \alpha \overline{(-\partial_x u(t,0) + \lambda (|u(t,0)|^{p-1} u(t,0)))}
	\int_{0}^{2 \pi} u(t,x) dx \bigg)\\
	&= \mathrm{Re}
	\bigg( \alpha \int_{0}^{2 \pi} u(t,x)
	\overline{(-\partial_x u(t,x) + \lambda (|u(t,x)|^{p-1} u(t,x)))} dx \bigg).
	\end{align*}
Then,
	\begin{align*}
	\frac{d}{dt} M(t)
	&= - \mathrm{Re} \thinspace \alpha
	\cdot \int_0^{2 \pi} 2 \mathrm{Re} ( u(t,x) \overline{\partial_x u(t,x) }) dx
	+ 2 \mathrm{Re} \thinspace \alpha \cdot \mathrm{Re} \thinspace \lambda
	\|u(t)\|_{L^{p+1}(\mathbb T)}^{p+1}\\
	&= - \mathrm{Re} \thinspace \alpha
	\cdot \int_0^{2\pi} \partial_x |u(t,x)|^2 dx
	+ 2 \mathrm{Re} \thinspace \alpha \cdot \mathrm{Re} \thinspace \lambda
	\|u(t)\|_{L^{p+1}(\mathbb T)}^{p+1}\\
	&=2 \mathrm{Re} \thinspace \alpha \cdot \mathrm{Re} \thinspace \lambda
	\|u(t)\|_{L^{p+1}(\mathbb T)}^{p+1}.
	\end{align*}
Since
	\[
	| M(t) |
	\leq | \alpha | \| u(t) \|_{L^1(\mathbb T)}^2
	\leq (2 \pi)^{\frac{2 p}{(p+1)}}
	| \alpha | \| u(t) \|_{L^{p+1}(\mathbb T)}^2,
	\]
we have
	\begin{align*}
	\frac{d}{dt} M(t)
	&\geq 2 (2 \pi)^{-p} | \alpha |^{-\frac{p+1}{2}}
	\mathrm{Re} \thinspace \alpha \cdot \mathrm{Re} \thinspace \lambda
	M(t)^{\frac{p+1}{2}}.
	\end{align*}
This implies
	\[
	M(t)
	\geq ( M(0)^{-\frac{p-1}{2}} - (p-1)(2 \pi)^{-p} | \alpha |^{-\frac{p+1}{2}}
	\mathrm{Re} \thinspace \alpha \cdot \mathrm{Re} \thinspace \lambda
	\thinspace t )^{- \frac{2}{p-1}}
	\]
and therefore
\begin{align*}
	T_0
	&\leq \inf \Big\{
	\frac{(2 \pi)^{p} | \alpha |^{\frac{p+1}{2}} }%
	{(p-1)\mathrm{Re} \thinspace \alpha \cdot \mathrm{Re} \thinspace \lambda}
	M(0)^{-\frac{p-1}{2}};\
	\alpha \in \mathbb C, \ \mathrm{Re} \thinspace \alpha \cdot \mathrm{Re} \lambda > 0
	\Big\}\\
	&\leq
	\frac{(2 \pi)^{p}}{(p-1) | \mathrm{Re} \thinspace \lambda|}
	\bigg| \int_0^{2 \pi} u_0(x) \int_0^x \overline{u_0(y)} dy \thinspace dx
	\bigg|^{-\frac{p-1}{2}}.
	\end{align*}



\end{document}